\documentclass[12pt]{article}
\usepackage{color}
\definecolor{darkblue}{rgb}{0.00,0.25,0.50}
\usepackage[colorlinks,filecolor=blue,citecolor=darkblue]{hyperref}

\usepackage{amsfonts,amssymb,amsmath,amsthm}
\usepackage{url}
\usepackage{enumerate}
\usepackage[ukrainian, russian, english]{babel}
\usepackage[cp1251]{inputenc}

\usepackage{cite}
\begin{document}
\selectlanguage{ukrainian} 
\thispagestyle{empty}

\title{}

\begin{center}
\textbf{\Large Оцінки знизу колмогоровських поперечників класів згорток з ядром Неймана}
\end{center}
\vskip0.5cm
\begin{center}
В.\,В.~Боденчук\\ \emph{\small Інститут математики НАН України, Київ}
\end{center}
\vskip0.5cm

\begin{abstract}
We obtain exact lower bounds for Kolmogorov $n$-widths in spaces $C$ and $L$ of classes of convolutions with Neumann kernel $N_{q,\beta}(t)=\sum\limits_{k=1}^{\infty}\dfrac{q^k}{k}\cos\left(kt-\dfrac{\beta\pi}{2}\right)$, ${q\in(0,1)}$, ${\beta\in\mathbb{R}}$, for all natural $n$ greater some number which depend only on $q$. The obtained estimates coincide with the best uniform approximations by trigonometric polynomials of mentioned classes. It made possible to obtain exact values for widths of these classes.

\vskip0.6cm

Встановлено точні оцінки знизу колмогоровських $n$-поперечників в просторах $C$ і $L$ класів згорток з ядром Неймана $N_{q,\beta}(t)=\sum\limits_{k=1}^{\infty}\dfrac{q^k}{k}\cos\left(kt-\dfrac{\beta\pi}{2}\right)$, ${q\in(0,1)}$, ${\beta\in\mathbb{R}}$, для усіх натуральних $n$, більших деякого номера, залежного лише від $q$. Отримані оцінки співпали з найкращими рівномірними наближеннями зазначених класів тригонометричними поліномами, що дозволило знайти точні значення поперечників цих класів.
\end{abstract}

\vskip1cm


\textbf{1. Постановка задачі.}
Через $L=L_1$ позначимо простір $2\pi$-періодичних сумовних функцій $f$ з нормою
$\|f\|_1=\int\limits_{-\pi}^{\pi}|f(t)|dt$, через
$L_\infty$ --- простір $2\pi$-періодичних вимірних і суттєво обмежених функцій з нормою
$\|f\|_\infty=\mathop{\rm ess\,sup}\limits_{t\in\mathbb{R}\ }|f(t)|,$
а через $C$ --- простір $2\pi$-періодичних неперервних функцій $f$, у якому норма задається рівністю
$\|f\|_C=\max\limits_{t\in\mathbb{R}}|f(t)|$.

Позначимо через $C_{\beta,p}^{\psi}$, $p=1, \infty$,  клас $2\pi$-періодичних функцій $f$, що зображуються у вигляді згортки
\begin{equation}\label{f}
f(x)=A+\dfrac{1}{\pi}\int\limits_{-\pi}^{\pi}\Psi_\beta(x-t)\varphi(t)dt=A+\left(\Psi_\beta\ast\varphi\right)(x),\; A\in\mathbb{R},
\end{equation}
де
\begin{equation}\label{Psi_beta}
\Psi_\beta(t)
	=
	\sum\limits_{k=1}^{\infty}
		\psi(k)\cos
			\left(
				kt
				-\dfrac{\beta\pi}{2}
			\right),
\,\psi(k)>0,
\, \sum\limits_{k=1}^{\infty}\psi(k)<\infty,\, \beta\in\mathbb{R},
\end{equation}
та
\begin{equation*}
\|\varphi\|_p\leqslant1,\, \varphi\perp1.
\end{equation*}
Функцію $\varphi$ в рівності \eqref{f} називають $(\psi,\beta)$-похідною функції $f$ і позначають через $f_\beta^\psi$.
Поняття $(\psi,\beta)$-похідної введене О.І. Степанцем (див., наприклад, \cite[\S~7--8]{Stepanets_2002_1}).

Якщо $\psi(k)=\dfrac{q^k}{k^r}$, $q\in(0,1)$, $r\in\mathbb{N}$, то класи $C_{\beta,p}^{\psi}$, $p=1, \infty$, позначають через $C_{\beta,p}^{q,r}$. При $r=1$ класи $C_{\beta,p}^{q,r}$ є класами згорток з ядром Неймана
\begin{equation*}
N_{q,\beta}(t)=\sum\limits_{k=1}^{\infty}\dfrac{q^k}{k}\cos\left(kt-\dfrac{\beta\pi}{2}\right), q\in(0,1), \beta\in\mathbb{R}.
\end{equation*}

Нехай $\mathfrak{N}$ --- центрально-симетрична підмножина банахового простору $X$. $m$-вимірним поперечником за Колмогоровим множини $\mathfrak{N}$ в просторі $X$ називається величина вигляду
\begin{equation}\label{d_m}
d_m(\mathfrak{N},X)=\inf_{F_m\subset X}\sup_{f\in \mathfrak{N}}\inf_{g\in F_m}\|f-g\|_X,
\end{equation}
де зовнішній $\inf$ береться по всіх $m$-вимірних лінійних підпросторах $F_m$ із $X$.

Розв’язанням задачі про знаходження оцінок колмогоровських поперечників функціональних класів займались багато математиків, із  результатами яких можна ознайомитись, наприклад, по монографіях \cite{Pinkus_1985,Tihomirov_1976,Kornejchuk_1987}.

У даній роботі розв’язується задача знаходження точних значень поперечників $d_{2n}(C_{\beta,\infty}^{q,1},C)$, $d_{2n-1}(C_{\beta,\infty}^{q,1},C)$ та $d_{2n-1}(C_{\beta,1}^{q,1},L)$  для усіх $q\in(0, 1)$, $\beta\in\mathbb{R}$ та натуральних $n$, більших деякого номера, залежного лише від $q$.

Для отримання оцінок поперечників зверху достатньо розглянути величину найкращого наближення класів $C_{\beta,p}^{\psi}$, де $p=1$ або $p=\infty$, в метриках просторів $X=L$ або $X=C$ відповідно підпростором $\mathcal{T}_{2n-1}$ тригонометричних поліномів $t_{n-1}$ порядку $n-1$, тобто величини вигляду
\begin{equation}\label{E_n}
E_n(C_{\beta,p}^{\psi})_X=\sup_{f\in C_{\beta,p}^{\psi}}\inf_{t_{n-1}\in \mathcal{T}_{2n-1}}\|f-t_{n-1}\|_X.
\end{equation}
Очевидно, що для величин вигляду \eqref{d_m} і \eqref{E_n} мають місце нерівності
\begin{equation}\label{d_m_E_n}
d_{2n-1}(C_{\beta,1}^{\psi},L)
	\leqslant 
	E_n(C_{\beta,1}^{\psi})_L.
\end{equation}
\begin{equation}\label{d_m_E_n2}
d_{2n-1}(C_{\beta,\infty}^{\psi},C)
	\leqslant 
	E_n(C_{\beta,\infty}^{\psi})_C.
\end{equation}

Як випливає з \cite{Serdyuk_1999}, для $\psi(k)=\dfrac{q^k}{k}$ при довільних $q\in(0,1)$, $\beta\in\mathbb{R}$ і $n\in\mathbb{N}$ мають місце рівності
\begin{equation*}
E_n(C_{\beta,\infty}^{q,1})_C
	=E_n(C_{\beta,1}^{q,1})_{L}
	=\|N_{q,\beta}\ast\varphi_n\|_C
	=
\end{equation*}
\begin{equation}\label{E_n_rivnosti}
	=\dfrac{4}{\pi}
		\left|
			\sum\limits_{\nu=0}^{\infty}
				\dfrac{q^{(2\nu+1)n}}{n(2\nu+1)^2}\sin
					\left(
						(2\nu+1)\theta_n\pi
						-\dfrac{\beta\pi}{2}
					\right)
		\right|,
\end{equation}
де
\begin{equation}\label{varp}
\varphi_n(t)=\textnormal{sign}\sin nt,
\end{equation}
а $\theta_n=\theta_n(q,\beta)$ --- єдиний на $[0,1)$ корінь рівняння
\begin{equation}\label{theta}
\sum\limits_{\nu=0}^{\infty}
	\dfrac
		{q^{2\nu n}}
		{2\nu+1}
	\cos\left(
		(2\nu+1)\theta_n\pi
		-
		\dfrac{\beta\pi}{2}
	\right)
	=0
\end{equation}
(при $\beta\in\mathbb{Z}$ рівності \eqref{E_n_rivnosti} випливають з \cite{Nikolskiy_1946, Nagy_1938}).

Згідно з \eqref{d_m_E_n}, \eqref{d_m_E_n2} і \eqref{E_n_rivnosti} задача про знаходження точних значень поперечників $d_{2n}(C_{\beta,\infty}^{q,1},C)$, $d_{2n-1}(C_{\beta,\infty}^{q,1},C)$ та $d_{2n-1}(C_{\beta,1}^{q,1},L)$ зводиться до встановлення справедливості оцінок знизу
\begin{equation}\label{dno1}
d_{2n}(C_{\beta,\infty}^{q,1}, C)\geqslant\|N_{q,\beta}\ast\varphi_n\|_C,
\end{equation}
\begin{equation}\label{dno2}
d_{2n-1}(C_{\beta,1}^{q,1}, L)\geqslant\|N_{q,\beta}\ast\varphi_n\|_C.
\end{equation}

При $q\in(0,1/7]$, $\beta\in\mathbb{Z}$ нерівності \eqref{dno1} та \eqref{dno2} випливають з результатів О.К.~Кушпеля \cite{Kushpel_1988} та \cite{Kushpel_1989}, а при $q\in(0, \,0{,}2]$, $\beta\in\mathbb{Z}$ і $q\in(0, \,0{,}193864]$, ${\beta\in\mathbb{R\setminus Z}}$ --- з роботи О.І.~Степанця та А.С.~Сердюка \cite{Serdyuk_1995}.
У вказаних випадках нерівності \eqref{dno1} та \eqref{dno2} доведено для усіх $n\in\mathbb{N}$.

\textbf{2. Основні результати.}
Для кожного фіксованого $q\in(0,1)$ позначимо через $n_{q}$ найменший з номерів $n\geqslant2$, для яких виконуються нерівності
\begin{equation}\label{umova_z}
\dfrac{q^{n}}{1-q^{2n}}\leqslant 
\min\{
\dfrac {2q^{\sqrt{n}}}{15n^2},
\dfrac{8}{3n^2} \left( \dfrac {2n-1}{7(n-1)^2}-\dfrac{\pi^2}{8n^2}\right)
\},
\end{equation}
\begin{equation}\label{umova_n_0}
\dfrac{24}{5(1-q)}q^{\sqrt{n}}+\dfrac{160}{63}\dfrac{2\sqrt{n}-1}{n(\sqrt{n}-1)}\; \dfrac{q}{(1-q)^2}\leqslant
\left(\dfrac{1}{2}+\dfrac{2q}{(1+q^2)(1-q)}\right)\left(\dfrac{1-q}{1+q}\right)^{\frac {4}{1-q^2}}.
\end{equation}

У прийнятих позначеннях має місце наступне твердження.

\textbf{Теорема 1.} \textit{Нехай $q\in(0,1)$.
Тоді для довільного $\beta\in \mathbb{R}$ і всіх номерів $n\geqslant n_{q}$ виконуються нерівності \eqref{dno1} і \eqref{dno2}}.

Як було відмічено вище, з роботи О.І.~Степанця та А.С.~Сердюка \cite{Serdyuk_1995} випливає, що при $q\in(0, \,0{,}2]$, $\beta\in\mathbb{Z}$ і $q\in(0, \,0{,}193864]$, ${\beta\in\mathbb{R\setminus Z}}$ нерівності \eqref{dno1} і \eqref{dno2} вірні для усіх $n\in\mathbb{N}$. 
Тоді, позначивши
\begin{equation*}
n_{q,\beta}\hspace{-3pt}=\hspace{-3pt}
\begin{cases}
1,
	&\hspace{-9pt}\text{якщо $q\in(0,\, 0{,}2]$ і $\beta\in\mathbb{Z}$ 
					 або $q\in(0,\,0{,}193864]$ і $\beta\in\mathbb{R\setminus Z}$,}
\\
n_{q},
	&\hspace{-9pt}\text{якщо $q\in(0{,}2,\,1)$ і $\beta\in\mathbb{Z}$ 
				   або $q\in(0{,}193864,\, 1)$ і $\beta\in\mathbb{R\setminus Z}$}
\end{cases}
\end{equation*}
одержимо істинність нерівностей \eqref{dno1} і \eqref{dno2} для усіх номерів $n\geqslant n_{q,\beta}$.
Об’єднавши нерівності \eqref{d_m_E_n}, \eqref{d_m_E_n2}, \eqref{dno1}, \eqref{dno2} з рівностями \eqref{E_n_rivnosti}
отримуємо наступне твердження.

\textbf{Теорема 2.}
\textit{Нехай $q\in(0,1)$ і $\beta\in\mathbb{R}$.
Тоді для довільного $\beta\in \mathbb{R}$ та усіх номерів $n\geqslant n_{q,\beta}$
мають місце рівності}
\begin{equation*}
d_{2n}(C_{\beta,\infty}^{q,1},C)
=d_{2n-1}(C_{\beta,\infty}^{q,1},C)
= d_{2n-1}(C_{\beta,1}^{q,1},L)
=
\end{equation*}
\begin{equation*}
=E_n(C_{\beta,\infty}^{q,1})_C
=E_n(C_{\beta,1}^{q,1})_L
=\|N_{q,\beta}\ast\varphi_n\|_C
=
\end{equation*}
\begin{equation}\label{dn}
=\dfrac{4}{\pi}
	\left|
		\sum\limits_{\nu=0}^{\infty}
			\dfrac{q^{(2\nu+1)n}}{n(2\nu+1)^2}
			\sin
				\left(
					(2\nu+1)\theta_n\pi
					-\dfrac{\beta\pi}{2}
				\right)
	\right|,
\end{equation}
\textit{де $\theta_n=\theta_n(q,\beta)$ --- єдиний на $[0,1)$ корінь рівняння \eqref{theta}}.

Теорема~2 дозволяє оцінити асимптотичну при ${n\to\infty}$ поведінку поперечників $d_{2n}(C_{\beta,\infty}^{q,1},C)$, $d_{2n-1}(C_{\beta,\infty}^{q,1},C)$ та $d_{2n-1}(C_{\beta,1}^{q,1},L)$.

\textbf{Теорема 3.}
\textit{Нехай $q\in(0,1)$ та $\beta\in \mathbb{R}$.
Тоді при $n\geqslant n_{q,\beta}$ }
\begin{equation*}
d_{2n}(C_{\beta,\infty}^{q,1},C)=d_{2n-1}(C_{\beta,\infty}^{q,1},C)= d_{2n-1}(C_{\beta,1}^{q,1},L)=E_n(C_{\beta,\infty}^{q,1})_C=
\end{equation*}
\begin{equation}\label{Th_3}
=E_n(C_{\beta,1}^{q,1})_L
=\dfrac{q^n}{n}
	\left(
		\dfrac{4}{\pi}
		+\gamma_n\frac{q^{2n}}{1-q^{2n}}
	\right),
\end{equation}
\textit{де $|\gamma_n|\leqslant\dfrac{16}{9\pi}$.}

Дійсно, знайдемо двосторонні оцінки правої частини формули \eqref{dn}. 
З \eqref{theta} отримуємо
\begin{equation*}
\left|
	\cos
	\left(
		\theta_n\pi
		-
		\dfrac{\beta\pi}{2}
	\right)
\right|
	=
		\left|
			\sum\limits_{\nu=1}^{\infty}
				\frac{q^{2\nu n}}{2\nu+1}
				\cos\left(
					(2\nu+1)\theta_n\pi
					-\dfrac{\beta\pi}{2}
				\right)
		\right|
	\leqslant
\end{equation*}
\begin{equation}\label{cos_y}
	\leqslant 
		\frac{1}{3} 
		\sum\limits_{\nu=1}^{\infty}
			q^{2\nu n}
		=\dfrac{q^{2n}}{3(1-q^{2n})}.
\end{equation}
З \eqref{cos_y} випливає, що
\begin{equation}\label{a_n}
0
	\leqslant
		 1-|\sin(\theta_n\pi-\dfrac{\beta\pi}{2})|
	\leqslant 
		|\cos(\theta_n\pi-\dfrac{\beta\pi}{2})|
	\leqslant 
		\dfrac{q^{2n}}{3(1-q^{2n})}.
\end{equation}

Оскільки,
\begin{equation*}
\left|
	\sum\limits_{\nu=1}^{\infty}
		\dfrac{q^{(2\nu+1)n}}{n(2\nu+1)^2}
		\sin
		\left(
			(2\nu+1)\theta_n\pi
			-\dfrac{\beta\pi}{2}
		\right)
\right|
	\leqslant
\end{equation*}
\begin{equation*}
	\leqslant
		\sum\limits_{\nu=1}^{\infty}
			\dfrac{q^{(2\nu+1) n}}{n(2\nu+1)^2}
	\leqslant
		\frac{1}{9n}\frac{q^{3n}}{1-q^{2n}},
\; n\in\mathbb{N},
\end{equation*}
то, враховуючи \eqref{a_n}, одержуємо для довільних $n\in\mathbb{N}$, $q\in(0,1)$ і $\beta\in\mathbb{R}$
\begin{equation*}
\left|
	\sum\limits_{\nu=0}^{\infty}
		\dfrac{q^{(2\nu+1)n}}{n(2\nu+1)^2}
		\sin
		\left(
			(2\nu+1)\theta_n\pi
			-\dfrac{\beta\pi}{2}
		\right)
\right|
	\geqslant
		\dfrac{q^n}{n}
		-\dfrac{q^n}{n}
		\left(
			1
			-|\sin(\theta_n\pi-\dfrac{\beta\pi}{2})|
		\right)
	-
\end{equation*}
\begin{equation*}
	-
	\left|
		\sum\limits_{\nu=1}^{\infty}
			\dfrac{q^{(2\nu+1)n}}{n(2\nu+1)^2}
			\sin
			\left(
				(2\nu+1)\theta_n\pi
				-\dfrac{\beta\pi}{2}
			\right)
	\right|
\geqslant
\end{equation*}
\begin{equation}\label{ots1}
\geqslant
	\dfrac{q^n}{n}
	\left(
		1
		-\frac{4}{9}\frac{q^{2n}}{1-q^{2n}}
	\right),
\end{equation}

\begin{equation*}
\left|
	\sum\limits_{\nu=0}^{\infty}
		\dfrac{q^{(2\nu+1)n}}{n(2\nu+1)^2}
		\sin
		\left(
			(2\nu+1)\theta_n\pi
			-\dfrac{\beta\pi}{2}
		\right)
\right|
	\leqslant
		\dfrac{q^n}{n}
		+\dfrac{q^n}{n}
		\left(
			1
			-|\sin(\theta_n\pi-\dfrac{\beta\pi}{2})|
		\right)
	+
\end{equation*}
\begin{equation*}
	+
	\left|
		\sum\limits_{\nu=1}^{\infty}
			\dfrac{q^{(2\nu+1)n}}{n(2\nu+1)^2}
			\sin
			\left(
				(2\nu+1)\theta_n\pi
				-\dfrac{\beta\pi}{2}
			\right)
	\right|
\leqslant
\end{equation*}
\begin{equation}\label{ots2}
\leqslant
	\dfrac{q^n}{n}
	\left(
		1
		+\frac{4}{9}\frac{q^{2n}}{1-q^{2n}}
	\right),
\end{equation}

З теореми~2 та оцінок \eqref{ots1} і \eqref{ots2} випливає, що при $n\geqslant n_{q,\beta}$ виконується \eqref{Th_3}.

Зазначимо, що оцінки \eqref{dno1} і \eqref{dno2} для довільних $q\in(0,1)$ неможливо встановити, використовуючи методи знаходження оцінок знизу для колмогоровських поперечників класів згорток із ядрами, що не збільшують осциляції розроблені А.~Пінкусом \cite{Pinkus_1979}, оскільки ядра Неймана $N_{q,\beta}(t)$ можуть збільшувати осциляцію.
Проілюструємо це на прикладі ядер Неймана $N_{q,0}(t)$ та $N_{q,1}(t)$ при $q=0{,}21$.

\textbf{Означення 1.}
 \textit{$2\pi$-періодичну функцію $K(\cdot)$ називають $\textnormal{\text{CVD}}_{2n}$-ядром (ядром, що не збільшує осциляції) і позначають $K\in \textnormal{\text{CVD}}_{2n}$, якщо для довільної функції $f\in C$ такої, що ${\nu(f)\leqslant 2n}$, виконується нерівність}
\begin{equation*}
\nu(K\ast f)\leqslant \nu(f),
\end{equation*}
 \textit{де $\nu(g)$--- число змін знаку функції $g\in C$ на $[0, 2\pi)$.}

Для встановлення факту, чи деяке ядро $\phi(x)$ є $\text{CVD}_{2n}$-ядром зручно користуватись наступним твердженням, що належить Мерхюберу, Шонбергу та Вільямсону \cite{Mairhuber} (див. також \cite[с.~67]{Pinkus_1985}).

\textit{Нехай $\phi(x)\in C$ та $\phi(x)$ має ранг не менший за
${2n+2}$, тобто існує розбиття $t_i$, $i=\overline{1,2n+2}$, проміжка $[0, 2\pi)$
таке, що $0\leqslant t_1<\dots<t_{2n+2}<2\pi$ і для якого
 $\dim(\textnormal{\text{span}}\{\phi(x-t_i)\}_{i=1}^{2n+2})=\linebreak =2n+2$. Тоді
  $\phi(x)\in \textnormal{\text{CVD}}_{2n}$ тоді і тільки тоді, коли
\begin{equation*}
D_{2l+1}(\mathbf{x},\mathbf{y})
	=\det(\varepsilon\phi(x_i-y_j))_{i,j=1}^{2l+1}
	\geqslant0,
\end{equation*}
\begin{equation*}
0\leqslant x_1<\dots<x_{2l+1}<2\pi,\;\;
0\leqslant y_1<\dots<y_{2l+1}<2\pi, \;\;
l=0, 1, \dots, n,
\end{equation*}
для деякого фіксованого   $\varepsilon=\pm1$.}

Як випливає із леми~1.3 роботи \cite{Kushpel_1985}, система функцій ${\{N_{q,\beta}(x-y_i)\}_{i=1}^{2n+2}}$ лінійно незалежна і, отже, ${\dim(\text{span}\{N_{q,\beta}(x-y_i)\}_{i=1}^{2n+2})=2n+2}$.
Тому, щоб довести, що ядра Неймана $N_{q,\beta}(t)$ при $q=0{,}21$ і $\beta=0$ або $\beta=1$ не є $\text{CVD}_{2n}$-ядрами ні при яких $n\in\mathbb{N}$, достатньо показати, що знайдуться вектори  $\mathbf{x}=(x_{1}, x_{2}, x_{3})$, $0\leqslant x_{1}<x_{2}<x_{3}<2\pi$, та $\mathbf{y}=(y_{1}, y_{2}, y_{3})$, $0\leqslant y_{1}<y_{2}<y_{3}<2\pi$, для яких детермінант $D_{3}(\mathbf{x},\mathbf{y})$ змінює знак. Виберемо вектори
 $\mathbf{x}^{(k)}=(x_{1}^{(k)}, x_{2}^{(k)}, x_{3}^{(k)})$ та $\mathbf{y}^{(k)}=(y_{1}^{(k)}, y_{2}^{(k)}, y_{3}^{(k)})$, $k=1, 2$, наступним чином:
\begin{equation*}
x_{1}^{(1)}=\frac{\pi}{18},\; x_{2}^{(1)}=\frac{\pi}{9},\; x_{3}^{(1)}=\frac{\pi}{6},\;
y_{1}^{(1)}=\frac{13\pi}{36},\; y_{2}^{(1)}=\frac{11\pi}{30},\; y_{3}^{(1)}=\frac{67\pi}{180},
\end{equation*}
\begin{equation*}
x_{1}^{(2)}=\frac{\pi}{18},\; x_{2}^{(2)}=\frac{\pi}{9},\; x_{3}^{(2)}=\frac{\pi}{6},\;
y_{1}^{(2)}=\frac{13\pi}{30},\; y_{2}^{(2)}=\frac{10\pi}{9},\; y_{3}^{(2)}=\frac{7\pi}{6}.
\end{equation*}
Обчислення показують, що для ядра $N_{q,0}$
\begin{equation*}
D_{3}(\mathbf{x}^{(1)},\mathbf{y}^{(1)})< -2{,}74\cdot 10^{-10},\;
D_{3}(\mathbf{x}^{(2)},\mathbf{y}^{(2)})> 1{,}09\cdot 10^{-6},
\end{equation*}
а для ядра $N_{q,1}$
\begin{equation*}
D_{3}(\mathbf{x}^{(1)},\mathbf{y}^{(1)})< -2{,}26\cdot 10^{-8},\;
D_{3}(\mathbf{x}^{(2)},\mathbf{y}^{(2)})> 2{,}09\cdot 10^{-6}.
\end{equation*}
Отже, для будь-яких $n\in\mathbb{N}$ при $q=0{,}21$\; $N_{q,0}(t)\not\in\text{CVD}_{2n}$ і $N_{q,1}(t)\not\in\text{CVD}_{2n}$.

\textbf{3. Означення і допоміжні твердження.} 
Для отримання оцінок \eqref{dno1} та \eqref{dno2} при довільних $\beta\in\mathbb{R}$ і $q\in(0, 1)$ в даній роботі буде використовуватись започаткований О.К.~Кушпелем \cite{Kushpel_1988} метод оцінки знизу поперечників класів згорток із твірними ядрами $\Psi_\beta$, що задовольняють так звану умову $C_{y,2n}$.
Наведемо необхідні означення та допоміжні твердження.

Нехай $\Delta_{2n}=\{0=x_0<x_1<\dots<x_{2n}=2\pi\}$, $x_k=k\pi/n$ --- розбиття проміжку $[0,2\pi]$. Розглянемо функцію
\begin{equation*}
\Psi_{\beta,1}(t)=(\Psi_{\beta}\ast B_{1})(t)=\sum\limits_{k=1}^{\infty}\frac{\psi(k)}{k}\cos\left(kt-\dfrac{(\beta+1)\pi}{2}\right),
\end{equation*}
де $B_{1}(t)=\sum\limits_{k=1}^{\infty}k^{-1}\sin kt$ --- ядро Бернуллі.
Через $S\Psi_{\beta,1}(\Delta_{2n})$ позначатимемо простір $SK$-сплайнів $S\Psi_{\beta,1}(\cdot)$ за розбиттям $\Delta_{2n}$, тобто множину функцій виду
\begin{equation}\label{SK}
S\Psi_{\beta,1}(\cdot)=\alpha_{0}+\sum\limits_{k=1}^{2n} \alpha_k \Psi_{\beta,1}(\cdot-x_k),\;\sum\limits_{k=1}^{2n} \alpha_k=0,
\end{equation}
\begin{equation*}
\alpha_k\in\mathbb{R},\; k= 0, 1, \dots, 2n.
\end{equation*}
Фундаментальним $SK$-сплайном називають функцію $\overline{S\Psi}_{\beta,1}(\cdot)=\overline{S\Psi}_{\beta,1}(y, \cdot)$ виду \eqref{SK}, що задовольняє співвідношення
\begin{equation*}
\overline{S\Psi}_{\beta,1}(y, y_k)=\delta_{0,k}=
\begin{cases}
0, & k=\overline{1,2n-1},    \\
1, & k=0,
\end{cases}
\end{equation*}
де $y_k=x_k+y$, $x_k=k\pi/n$, $y\in[0,\dfrac{\pi}{n})$. 
Оскільки серед $(\psi, \beta)$-похідних будь-якого сплайна виду \eqref{SK} існує функція, яка є сталою на кожному інтервалі $(x_k, x_{k+1})$, то надалі саме таку функцію будемо розуміти під записом $(\overline{S\Psi}_{\beta,1}(\cdot))_\beta^\psi$.

\textbf{Означення 2.} \textit{Будемо казати, що для деякого дійсного числа $y$ і розбиття $\Delta_{2n}$ ядро $\Psi_\beta(\cdot)$ вигляду \eqref{Psi_beta} задовольняє умову $C_{y,2n}$ (і записувати $\Psi_\beta\in C_{y,2n}$), якщо для цього ядра існує єдиний фундаментальний сплайн $\overline{S\Psi}_{\beta,1}(y,\cdot)$ і для нього виконуються рівності}
\begin{equation*}
\textrm{sign}(\overline{S\Psi}_{\beta,1}(y,t_k))_\beta^\psi=(-1)^k\varepsilon e_k,\, k=\overline{0,2n-1},
\end{equation*}
\textit{де $t_k=(x_{k}+x_{k+1})/2,$ $e_k$ дорівнює або 0, або 1, а $\varepsilon$ приймає значення $\pm1$ і не залежить від $k$.}

Наступна теорема дозволяє знаходити оцінки знизу колмогоровських поперечників класів згорток, породжених ядрами, що задовольняють умову $C_{y,2n}$.

\textbf{Теорема 4 (О.К.~Кушпель \cite{Kushpel_1988,Kushpel_1989}).} \textit{Нехай при деякому $n\in\mathbb{N}$ функція $\Psi_{\beta}$ вигляду \eqref{Psi_beta}, що породжує класи $C_{\beta,p}^\psi$, $p=1,\infty$, задовольняє умову $C_{y,2n}$, коли $y$ --- точка, в якій функція $|(\Psi_\beta\ast\varphi_n)(t)|$, $\varphi_n(t)=\textnormal{sign}\sin nt$, приймає максимальне значення. Тоді }
\begin{equation*}
d_{2n}(C_{\beta,\infty}^\psi, C)\geqslant\|\Psi_\beta\ast\varphi_n\|_C,
\end{equation*}
\begin{equation*}
d_{2n-1}(C_{\beta,1}^\psi, L)\geqslant\|\Psi_\beta\ast\varphi_n\|_C.
\end{equation*}

У роботах \cite{Kushpel_1988,Serdyuk_1998,Serdyuk_1995,Serdyuk_1999} були встановлені достатні умови включення $\Psi_\beta\in C_{y,2n}$ для ядер виду \eqref{Psi_beta}. Це дозволило авторам зазначених робіт застосувати теорему 4 і одержати в ряді нових випадків точні оцінки поперечників $d_m(C_{\beta,\infty}^\psi,C)$ та $d_m(C_{\beta,1}^\psi,L)$.

\textbf{Лема 1.} \textit{Нехай $\beta\in \mathbb{R}$, $\sum\limits_{k=1}^{\infty}\psi(k)<\infty$ і  
\begin{equation}\label{lambda_not0}
|\lambda_{l}(y)|\not=0,\;l=\overline{1,n},
\end{equation}
де
\begin{equation}\label{lambd}
\lambda_{l}(y)=\dfrac1n \sum_{\nu=1}^{2n}e^{il\nu\pi/n}\Psi_{\beta,1}(y-\dfrac{\nu\pi}{n}).
\end{equation}
Тоді для довільного ${t\in(\dfrac{(k-1)\pi}{n},\dfrac{k\pi}{n})}$, $k=\overline{1,2n}$, виконується рівність}
\begin{equation*}
(\overline{S\Psi}_{\beta,1}(y,t))_\beta^\psi=
(-1)^{k+1}\frac{\pi}{4n\psi(n)}
\times
\end{equation*}
\begin{equation}\label{SP_Psi}
\times
\Bigg(
	\bigg(
	\dfrac{1}{2}
	+2\dfrac{\psi(n)}{n}\sum_{j=1}^{n-1}\dfrac{\cos j(t_k-y)}{|\lambda_{n-j}(y)|\cos\dfrac{j\pi}{2n}}
	\bigg)
	\mathop{\text{sign}} \sin(ny-\frac{\beta\pi}{2})
	+\gamma_1(y)
	+\gamma_2(y)
\Bigg),
\end{equation}
\textit{в якій $t_k=\dfrac{k \pi}{n}-\dfrac{\pi}{2n}$, а } 
\begin{equation*}
\gamma_1(y)
	=
		\gamma_1(\psi,\beta,k,y)
	=
\end{equation*}
\begin{equation}\label{gamma_2}
	=
		\frac{\psi(n)}{n}
		\left(
			\dfrac{z_{0}(y)}{|\lambda_{n}(y)|^2}+
			2\sum_{j=1}^{n-1}
				\dfrac{z_{j}(y)}{|\lambda_{n-j}(y)|^2\cos\dfrac{j\pi}{2n}}
		\right),
\end{equation}
\begin{equation}\label{gamma_3}
\gamma_2(y)
	=
		\gamma_2(\psi,\beta,y)
	=
		-
		\dfrac
			{R_0(y)\dfrac{n}{\psi(n)}}
			{2(2+R_0(y)\dfrac{n}{\psi(n)})}
		\mathop{\text{sign}} \sin(ny-\frac{\beta\pi}{2}),
\end{equation}
\begin{equation*}
z_{j}(y)
	=
		z_{j}(\psi,\beta,k,y)
	=
		|r_{j}(y)|\cos(j(t_k-y)+\arg(r_{j}(y)))
		-
\end{equation*}
\begin{equation}\label{z_nj}
		-
		R_{j}(y)\cos(j(t_k-y))
		\mathop{\text{sign}} \sin(ny-\frac{\beta\pi}{2}),
		\;j=\overline{0,n-1},
\end{equation}
\begin{equation}\label{|lambda_n-j|}
R_{j}(y)
	=
		R_{j}(\psi,\beta,y)
	=
		|\lambda_{n-j}(y)|
		-
		\dfrac{\psi(n-j)}{n-j}
		-
		\dfrac{\psi(n+j)}{n+j},
		\;j=\overline{0,n-1},
\end{equation}
\begin{align}
\label{r_n-j}
&r_{j}(y)
	=
		\sum_{\nu=1}^3
			r_{j}^{(\nu)}(y),
		\; j=\overline{0,n-1},
\\
\nonumber
&r_{j}^{(1)}(y)
	=
		r_{j}^{(1)}(\psi,\beta,y)
	=
		\dfrac
			{\psi(3n-j)e^{i(3ny-\frac{(\beta+1)\pi}{2})}}
			{3n-j}
		+
\\
\nonumber
&\phantom{r_{j}^{(1)}(y)=}
	+
	\sum\limits_{m=2}^{\infty}
		\left(
			\dfrac
				{\psi((2m+1)n-j)e^{i((2m+1)ny-\frac{(\beta+1)\pi}{2})}}
				{(2m+1)n-j}
			+
		\right.
\\
\label{r_nj_1}
&\phantom{r_{j}^{(1)}(y)=}
		\left.
			+
			\dfrac
				{\psi((2m-1)n+j)e^{-i((2m-1)ny-\frac{(\beta+1)\pi}{2})}}
				{(2m-1)n+j}
		\right),
\\
\label{r_nj_2}
&r_{j}^{(2)}(y)
	=
		r_{j}^{(2)}(\psi,\beta,y)
	=
		i
		\left(
			\dfrac{\psi(n+j)}{n+j}
			-
			\dfrac{\psi(n-j)}{n-j}
		\right)
		\cos(ny-\frac{\beta\pi}{2}),
\\
\nonumber
&r_{j}^{(3)}(y)
	=
		r_{j}^{(3)}(\psi,\beta,y)
	=
		\left(
			\dfrac{\psi(n-j)}{n-j}
			+
			\dfrac{\psi(n+j)}{n+j}
		\right)
		\times	
\\
\label{r_nj_3}
&\phantom{r_{j}^{(3)}(y)=}
	\times		
		(|\sin(ny-\dfrac{\beta\pi}{2})|-1)
		\mathop{\text{sign}} \sin(ny-\dfrac{\beta \pi}{2}).
\end{align}

\textbf{Доведення.} Будемо виходити із отриманого у роботі \cite{Serdyuk_1995} зображення функції $(\overline{S\Psi}_{\beta,1}(y,t))_\beta^\psi$, згідно з яким за умови $|\lambda_j(y)|\not=0$, $j=\overline{1, n}$, для довільного $t\in(x_{k-1},x_k)$ виконується рівність
\begin{equation}\label{SPsi_v0}
(\overline{S\Psi}_{\beta,1}(y,t))_\beta^\psi=\frac{\pi}{4n^2}\left(2\sum_{j=1}^{n-1}\dfrac{\sin jt_k\cdot\rho_j(y)-\cos jt_k\cdot\sigma_j(y)}{|\lambda_j(y)|^2\sin\dfrac{j\pi}{2n}}+\dfrac{(-1)^{k+1}\rho_n(y)}{|\lambda_n(y)|^2}\right),
\end{equation}
де
\begin{equation*}
\lambda_j(\cdot)=\dfrac1n \sum_{\nu=1}^{2n}e^{ij\nu\pi/n}\Psi_{\beta,1}(\cdot-\dfrac{\nu\pi}{n}),
\end{equation*}
$i$ --- уявна одиниця, $\rho_j(\cdot)=\mathop{\text{Re}}(\lambda_j(\cdot))$, $\sigma_j(\cdot)=\mathop{\text{Im}}(\lambda_j(\cdot))$, $t_k=\dfrac{k \pi}{n}-\dfrac{\pi}{2n}$.

Змінивши порядок підсумовування доданків у сумі в правій частині рівності \eqref{SPsi_v0}, маємо
\begin{equation*}
\sum_{j=1}^{n-1}\dfrac{\sin jt_k\cdot\rho_j(y)-\cos jt_k\cdot\sigma_j(y)}{|\lambda_j(y)|^2\sin\dfrac{j\pi}{2n}}=
\end{equation*}
\begin{equation*}
=\sum_{j=1}^{n-1}\dfrac{\sin (n-j)t_k\cdot\rho_{n-j}(y)-\cos (n-j)t_k\cdot\sigma_{n-j}(y)}{|\lambda_{n-j}(y)|^2\sin\dfrac{(n-j)\pi}{2n}}=
\end{equation*}
\begin{equation}\label{S_Psi}
=(-1)^{k+1}\sum_{j=1}^{n-1}\dfrac{\cos jt_k\cdot\rho_{n-j}(y)-\sin jt_k\cdot\sigma_{n-j}(y)}{|\lambda_{n-j}(y)|^2\cos\dfrac{j\pi}{2n}}.
\end{equation}

З урахуванням \eqref{SPsi_v0}  і \eqref{S_Psi} для фундаментального $SK$-сплайна $\overline{S\Psi}_{\beta,1}(y,t)$, 
 за умови $|\lambda_j(y)|\not=0$, $j=\overline{1, n}$, одержуємо зображення
\begin{equation*}
(\overline{S\Psi}_{\beta,1}(y,t))_\beta^\psi=
\end{equation*}
\begin{equation}\label{SP_v0}
=\frac{(-1)^{k+1}\pi}{4n^2}\left(2\sum_{j=1}^{n-1}\dfrac{\cos jt_k\cdot\rho_{n-j}(y)-\sin jt_k\cdot\sigma_{n-j}(y)}{|\lambda_{n-j}(y)|^2\cos\dfrac{j\pi}{2n}}+\dfrac{\rho_n(y)}{|\lambda_n(y)|^2}\right).      
\end{equation}

Покажемо, що величини $\lambda_{n-j}(y)$ виду \eqref{lambd} при $j=\overline{0,n-1}$ можна виразити наступним чином:
\begin{equation}\label{lambda_n-j}
\lambda_{n-j}(y)=
e^{-ijy}
\left(
	\left(
		\dfrac{\psi(n-j)}{n-j}+\dfrac{\psi(n+j)}{n+j}
	\right)
	\mathop{\text{sign}} \sin(ny-\dfrac{\beta \pi}{2})
	+r_{j}(y)
\right),
\end{equation}
де величини $r_{j}(y)$ задаються рівностями \eqref{r_n-j}.

Перепишемо ядро $\Psi_{\beta,1}$ у комплексній формі
\begin{equation*}
\Psi_{\beta,1}(t)
=(\Psi_{\beta}\ast B_1)(t)
=\sum_{k=1}^\infty\frac{\psi(k)}{k}\cos (kt-\frac{(\beta+1)\pi}{2})
= \dfrac 12 \sideset{}{'}\sum\limits_{k=-\infty}^{\infty}c_ke^{ikt},
\end{equation*}
де
\begin{equation}\label{c_k}
c_k=\dfrac{\psi(k)}{k}e^{-i\frac{(\beta+1)\pi}{2}},\;c_{-k}=\dfrac{\psi(k)}{k}e^{i\frac{(\beta+1)\pi}{2}},\;k\in\mathbb{N},
\end{equation}
а штрих біля знака суми означає, що при підсумовуванні відсутній доданок з нульовим номером.

Підставивши у \eqref{lambd} замість ядра $\Psi_{\beta,1}$ його розклад у комплексний ряд Фур’є, одержимо
\begin{equation*}
\lambda_l(y)=\dfrac1n \sum_{\nu=1}^{2n}e^{il\nu\pi/n}\dfrac 12 \sideset{}{'}\sum\limits_{k=-\infty}^{\infty}c_ke^{ik(y-\nu\pi/n)} =
\end{equation*}
\begin{equation*}
=\dfrac{1}{2n} \sum_{\nu=1}^{2n}\sideset{}{'}\sum\limits_{k=-\infty}^{\infty}c_k e^{i(ky+(l-k)\nu\pi/n)} =
\end{equation*}
\begin{equation}\label{eq:lambda_dop1}
=\dfrac{1}{2n} \sideset{}{'}\sum\limits_{k=-\infty}^{\infty}c_k e^{iky} \sum_{\nu=1}^{2n}e^{i((l-k)\nu\pi/n)}.
\end{equation}

Неважко переконатись, що
\begin{equation}\label{eq:lambda_dop2}
\sum_{\nu=1}^{2n}e^{i((l-k)\nu\pi/n)}=\begin{cases}
0,& \text{якщо } k\not=l-2mn, m\in\mathbb{Z};\\
2n,& \text{якщо } k=l-2mn, m\in\mathbb{Z}.
\end{cases}
\end{equation}

З \eqref{eq:lambda_dop1} та \eqref{eq:lambda_dop2} при $l=\overline{1,n}$ випливає наступне представлення:
\begin{equation*}
\lambda_l(y)=\sum\limits_{m=-\infty}^{+\infty}c_{l-2mn}e^{i(l-2mn)y}
=\sum\limits_{m=-\infty}^{+\infty}c_{2mn+l}e^{i(2mn+l)y}.
\end{equation*}

Звідси при $l=n-j\,$,   $j=\overline{0,n-1}$, отримуємо
\begin{equation*}
\lambda_{n-j}(y)=\sum\limits_{m=-\infty}^{+\infty}c_{(2m+1)n-j}e^{i((2m+1)n-j)y}=
\end{equation*}
\begin{equation}\label{lambda_n-j_y0}
=e^{-ijy}(c_{n-j}e^{iny}+c_{-(n+j)}e^{-iny}+r_{j}^{(1)}(y)).
\end{equation}

З урахуванням \eqref{c_k} перетворимо перші два доданки в \eqref{lambda_n-j_y0} наступним чином:
\begin{equation*}
c_{n-j}e^{iny}+c_{-(n+j)}e^{-iny}=
\end{equation*}
\begin{equation*}
=\dfrac{\psi(n-j)}{n-j}e^{i(ny-\frac{(\beta+1)\pi}{2})}
+\dfrac{\psi(n+j)}{n+j}e^{-i(ny-\frac{(\beta+1)\pi}{2})}=
\end{equation*}
\begin{equation*}
=\left(\dfrac{\psi(n-j)}{n-j}+\dfrac{\psi(n+j)}{n+j}\right)\cos(ny-\frac{(\beta+1)\pi}{2})+
\end{equation*}
\begin{equation*}
+i\left(\dfrac{\psi(n-j)}{n-j}-\dfrac{\psi(n+j)}{n+j}\right)\sin(ny-\frac{(\beta+1)\pi}{2})=
\end{equation*}
\begin{equation}\label{c_n-j-}
=\left(\dfrac{\psi(n-j)}{n-j}+\dfrac{\psi(n+j)}{n+j}\right)\sin(ny-\frac{\beta\pi}{2})+r_{j}^{(2)}(y).
\end{equation}

Записавши $\sin(ny-\dfrac{\beta\pi}{2})$ у вигляді
\begin{equation*}\label{sin_ny}
\sin\left(ny-\dfrac{\beta\pi}{2}\right)=
|\sin(ny-\dfrac{\beta\pi}{2})|\,
\mathop{\text{sign}} \sin(ny-\dfrac{\beta \pi}{2}),
\end{equation*}
з \eqref{c_n-j-} маємо
\begin{equation*}
c_{n-j}e^{iny}+c_{-(n+j)}e^{-iny}=
\left(
	\dfrac{\psi(n-j)}{n-j}+\dfrac{\psi(n+j)}{n+j}
\right)
\mathop{\text{sign}} \sin(ny-\frac{\beta\pi}{2})+
\end{equation*}
\begin{equation*}
+\left(
	\dfrac{\psi(n-j)}{n-j}+\dfrac{\psi(n+j)}{n+j}
\right)
(|\sin(ny-\dfrac{\beta\pi}{2})|-1)
\mathop{\text{sign}} \sin(ny-\dfrac{\beta \pi}{2})
+r_{j}^{(2)}(y)=
\end{equation*}
\begin{equation}\label{c_n-j}
=\left(
	\dfrac{\psi(n-j)}{n-j}
	+\dfrac{\psi(n+j)}{n+j}
\right)
\mathop{\text{sign}} \sin(ny-\frac{\beta\pi}{2})
+r_{j}^{(2)}(y)
+r_{j}^{(3)}(y).
\end{equation}

Рівності \eqref{lambda_n-j_y0} та \eqref{c_n-j} доводять формулу \eqref{lambda_n-j}.

Перетворимо чисельник кожного доданка в правій частині рівності (\ref{SP_v0}). Для цього, з урахуванням (\ref{lambda_n-j}), запишемо
\begin{equation*}
\rho_{n-j}(y)= \mathop{\text{Re}}(\lambda_{n-j}(y))=
\end{equation*}
\begin{equation}\label{rho_n-j}
=\left(\dfrac{\psi(n-j)}{n-j}+\dfrac{\psi(n+j)}{n+j}\right)\cos jy \;
\mathop{\text{sign}} \sin(ny-\frac{\beta\pi}{2})
+\mathop{\text{Re}}(e^{-ijy}r_{j}(y));
\end{equation}
\begin{equation*}
\sigma_{n-j}(y)=\mathop{\text{Im}}(\lambda_{n-j}(y))=
\end{equation*}
\begin{equation}\label{sigma_n-j}
=-\left(\dfrac{\psi(n-j)}{n-j}+\dfrac{\psi(n+j)}{n+j}\right)\sin jy \;
\mathop{\text{sign}} \sin(ny-\frac{\beta\pi}{2})
+\mathop{\text{Im}}(e^{-ijy}r_{j}(y)).
\end{equation}

Застосовуючи \eqref{rho_n-j} та \eqref{sigma_n-j}, отримуємо
\begin{equation*}
\cos jt_k\cdot\rho_{n-j}(y)
-\sin jt_k\cdot\sigma_{n-j}(y)=
\end{equation*}
\begin{equation*}
=
\left(
	\dfrac{\psi(n-j)}{n-j}
	+\dfrac{\psi(n+j)}{n+j}
\right)
\cos(j(t_k-y))
\,\mathop{\text{sign}} \sin(ny-\frac{\beta\pi}{2})+
\end{equation*}
\begin{equation*}
+\cos jt_k\cdot\mathop{\text{Re}} (e^{-ijy} r_{j}(y))
-\sin jt_k\cdot\mathop{\text{Im}} (e^{-ijy} r_{j}(y))=
\end{equation*}
\begin{equation*}
=
\left(
	\dfrac{\psi(n-j)}{n-j}
	+\dfrac{\psi(n+j)}{n+j}
	+R_{j}(y)
\right)
\cos(j(t_k-y)) 
\,\mathop{\text{sign}} \sin(ny-\frac{\beta\pi}{2})
+
\end{equation*}
\begin{equation}\label{cos+sin}
+z_{j}(y)=
|\lambda_{n-j}(y)|\cos(j(t_k-y)) 
\,\mathop{\text{sign}} \sin(ny-\frac{\beta\pi}{2})
+z_{j}(y),
\end{equation}
де
\begin{equation*}
z_{j}(y)=
\cos jt_k\cdot\mathop{\text{Re}}(e^{-ijy} r_{j}(y))
-\sin jt_k\cdot\mathop{\text{Im}} (e^{-ijy} r_{j}(y))-
\end{equation*}
\begin{equation*}
-R_{j}(y)\cos(j(t_k-y))
\,\mathop{\text{sign}} \sin(ny-\frac{\beta\pi}{2}),
\end{equation*}
а $R_{j}(y)$ означені у \eqref{|lambda_n-j|}.

В силу очевидної рівності
\begin{equation*}
e^{-ijy} r_{j}(y)=|r_{j}(y)|(\cos (\arg(r_{j}(y))-jy)+i\sin (\arg(r_{j}(y))-jy))
\end{equation*}
величину $z_{j}(y)$ можна зобразити у вигляді \eqref{z_nj}.

При $j=0$ формула \eqref{lambda_n-j} перетворюється в наступну рівність:
\begin{equation}\label{lambda_n}
\lambda_{n}(y)=
2\dfrac{\psi(n)}{n}
\mathop{\text{sign}} \sin(ny-\frac{\beta\pi}{2})
+r_{0}(y),
\end{equation}
де $r_{0}(y)$ визначається формулою \eqref{r_n-j}, у якій
\begin{equation}\label{r_n1}
r_{0}^{(1)}(y)=
2\sum\limits_{m=2}^{\infty}\dfrac{\psi((2m-1)n)}{(2m-1)n}\cos((2m-1)ny-\frac{(\beta+1)\pi}{2}),
\end{equation}
\begin{equation}\label{r_n2}
r_{0}^{(2)}(y)=0,
\end{equation}
\begin{equation}\label{r_n3}
r_{0}^{(3)}(y)=2\dfrac{\psi(n)}{n} (|\sin(ny-\dfrac{\beta\pi}{2})|-1)
\mathop{\text{sign}} \sin(ny-\frac{\beta\pi}{2}).
\end{equation}

З \eqref{lambda_n}--\eqref{r_n3} випливає, що $\sigma_{n}(y)=0$ і тому
\begin{equation*}
\rho_{n}(y)=\lambda_{n}(y)=
2\dfrac{\psi(n)}{n}
\mathop{\text{sign}} \sin(ny-\frac{\beta\pi}{2})
+r_{0}(y).  
\end{equation*}

Звідси, враховуючи \eqref{z_nj} та \eqref{|lambda_n-j|}, можна записати
\begin{equation*}
\rho_{n}(y)=
\left(
	2\dfrac{\psi(n)}{n}
	+R_0(y)
\right)
\mathop{\text{sign}} \sin(ny-\frac{\beta\pi}{2})
+z_{0}(y)=
\end{equation*}
\begin{equation}\label{rho_n}
=|\lambda_{n}(y)|
\mathop{\text{sign}} \sin(ny-\frac{\beta\pi}{2})
+z_{0}(y),
\end{equation}
де
\begin{equation*}
z_0(y)=r_0(y)- R_0(y)
\mathop{\text{sign}} \sin(ny-\frac{\beta\pi}{2}).
\end{equation*}

Із зображення (\ref{SP_v0}) і рівностей (\ref{cos+sin}) та \eqref{rho_n} отримуємо
\begin{equation*}
(\overline{S\Psi}_{\beta,1}(y,t))_\beta^\psi=
\end{equation*}
\begin{equation*}
=\frac{(-1)^{k+1}\pi}{4n\psi(n)}
\left(
	\mathop{\text{sign}} \sin(ny-\frac{\beta\pi}{2})
	\left(
		2\frac{\psi(n)}{n}\sum_{j=1}^{n-1}\dfrac{\cos j(t_k-y)}{|\lambda_{n-j}(y)|\cos\dfrac{j\pi}{2n}}
		+\dfrac{\psi(n)}{n|\lambda_n(y)|}
	\right)
	+
\right.
\end{equation*}
\begin{equation*}
\left.
	+2\frac{\psi(n)}{n}\sum_{j=1}^{n-1}\dfrac{z_{j}(y)}{|\lambda_{n-j}(y)|^2\cos\dfrac{j\pi}{2n}}
	+\dfrac{\psi(n)z_{0}(y)}{n|\lambda_{n}(y)|^2}
\right)=
\end{equation*}
\begin{equation*}
=\frac{(-1)^{k+1}\pi}{4n\psi(n)}
\Bigg(
	\mathop{\text{sign}} \sin(ny-\frac{\beta\pi}{2})
\Bigg.
\times
\end{equation*}
\begin{equation}\label{Phi_x}
\times
\Bigg.
	\Bigg(
		2\frac{\psi(n)}{n}\sum_{j=1}^{[\sqrt{n}]}\dfrac{\cos j(t_k-y)}{|\lambda_{n-j}(y)|\cos\dfrac{j\pi}{2n}}
		+\dfrac{\psi(n)}{n|\lambda_n(y)|}
	\Bigg)
	+\gamma_1(y)
\Bigg).
\end{equation}

В силу (\ref{|lambda_n-j|})
\begin{equation*}
\dfrac{\psi(n)\mathop{\text{sign}} \sin(ny-\dfrac{\beta\pi}{2})}{n|\lambda_n(y)|}=
\dfrac{\mathop{\text{sign}} \sin(ny-\dfrac{\beta\pi}{2})}{2+R_0(y)\dfrac{n}{\psi(n)}}=
\end{equation*}
\begin{equation*}
=\Bigg(
	\dfrac{1}{2}
	-\dfrac{R_0(y)\dfrac{n}{\psi(n)}}{2(2+R_0(y)\dfrac{n}{\psi(n)})}
\Bigg)
\mathop{\text{sign}} \sin(ny-\frac{\beta\pi}{2})=
\end{equation*}
\begin{equation}\label{s0}
=\dfrac{1}{2}\mathop{\text{sign}} \sin(ny-\frac{\beta\pi}{2})+\gamma_2(y).
\end{equation}

Із \eqref{Phi_x} та \eqref{s0} отримуємо \eqref{SP_Psi}.
Лему доведено.

\textbf{Лема 2.} \textit{Нехай коефіцієнти $\psi(k)>0$ ядра $\Psi_\beta$ задовільняють умову}
\begin{equation*}
\lim_{k\to\infty}\dfrac{\psi(k+1)}{\psi(k)}=q,
	q\in(0,1),
\end{equation*}
\textit{та $\beta\in \mathbb{R}$.
Тоді при виконанні умови \eqref{lambda_not0} для довільного ${t\in(\dfrac{(k-1)\pi}{n},\dfrac{k\pi}{n})}$, $k=\overline{1,2n}$, справедлива рівність}
\begin{equation}\label{SP_Phi_}
(\overline{S\Psi}_{\beta,1}(y,t))_\beta^\psi=
(-1)^{k+1}\frac{\pi}{4n\psi(n)}\;
\bigg(
	\mathcal{P}_q(t_k-y)
	\mathop{\text{sign}} \sin(ny-\frac{\beta\pi}{2})
	+\sum_{m=1}^5\gamma_m(y)
\bigg),
\end{equation}
\textit{в якій $t_k=\dfrac{k \pi}{n}-\dfrac{\pi}{2n}$, $\mathcal{P}_q(t)$ --- ядро аналітично продовжуваних в смугу функцій,}
\begin{equation*}
\mathcal{P}_q(t)= \dfrac{1}{2}+2\sum_{j=1}^{\infty}\dfrac{\cos jt}{q^j+q^{-j}},
\end{equation*}
\textit{величини $\gamma_1(y)$ та $\gamma_2(y)$ задані рівностями \eqref{gamma_2} і \eqref{gamma_3} відповідно, а}
\begin{equation}\label{gamma_1}
\gamma_3(y)
	=
		\gamma_3(\psi,\beta,k,y)
	=
		2\sum_{j=[\sqrt{n}]+1}^{n-1}
			\dfrac
				{\cos j(t_k-y)}
				{\dfrac{n}{\psi(n)}|\lambda_{n-j}(y)|\cos\dfrac{j\pi}{2n}}
		\mathop{\text{sign}} \sin(ny-\frac{\beta\pi}{2}),
\end{equation}
\begin{equation}\label{gamma_4}
\gamma_4(y)
	=
		\gamma_4(\psi,\beta,k,y)
	=
		-2\sum_{j=1}^{[\sqrt{n}]}
			\dfrac
				{\delta_{j}(y)\cos j(t_k-y)}
				{\dfrac{n}{\psi(n)}|\lambda_{n-j}(y)|\cos\dfrac{j\pi}{2n}}
			\mathop{\text{sign}} \sin(ny-\frac{\beta\pi}{2}),
\end{equation}
\begin{equation}\label{gamma_5}
\gamma_5(y)
	=
		\gamma_5(q,\beta,k,y)
	=
		-2\sum\limits_{j=[\sqrt{n}]+1}^{\infty}
			\dfrac
				{\cos j(t_k-y)}
				{q^j+q^{-j}}
			\mathop{\text{sign}} \sin(ny-\frac{\beta\pi}{2}),
\end{equation}
\begin{equation}\label{delta_0}
\delta_{j}(y)
	=
		\delta_{j}(\psi,y)
	=
	\dfrac
		{n|\lambda_{n-j}(y)|\cos\dfrac{j\pi}{2n}}
		{(q^{-j}+q^{j})\psi(n)}
	-1,
	\;j=\overline{1,[\sqrt{n}]},
\end{equation}
\textit{$[a]$ --- ціла частина числа $a$.}

\textbf{Доведення.} Згідно з \eqref{gamma_1}
\begin{equation*}
2\frac{\psi(n)}{n}\sum_{j=1}^{n-1}
	\dfrac{\cos j(t_k-y)}{|\lambda_{n-j}(y)|\cos\dfrac{j\pi}{2n}}
	\mathop{\text{sign}} \sin(ny-\frac{\beta\pi}{2})
	=
\end{equation*}
\begin{equation}\label{summa_gamma}
=2\frac{\psi(n)}{n}\sum_{j=1}^{[\sqrt{n}]}\dfrac{\cos j(t_k-y)}{|\lambda_{n-j}(y)|\cos\dfrac{j\pi}{2n}}
\mathop{\text{sign}} \sin(ny-\frac{\beta\pi}{2})
+\gamma_3(y).
\end{equation}

Далі в силу формул \eqref{gamma_4}--\eqref{delta_0} можна записати рівності
\begin{equation*}
2\dfrac{\psi(n)}{n}\sum_{j=1}^{[\sqrt{n}]}\dfrac{\cos j(t_k-y)}{|\lambda_{n-j}(y)|\cos\dfrac{j\pi}{2n}}
\mathop{\text{sign}} \sin(ny-\frac{\beta\pi}{2})=
\end{equation*}
\begin{equation*}
=
2\sum_{j=1}^{[\sqrt{n}]}\dfrac{\cos j(t_k-y)}{(q^j+q^{-j})(1+\delta_{j}(y))}
\mathop{\text{sign}} \sin(ny-\frac{\beta\pi}{2})=
\end{equation*}
\begin{equation*}
=
\left(
	2\sum_{j=1}^{[\sqrt{n}]}
		\dfrac{\cos j(t_k-y)}{q^j+q^{-j}}
	-2\sum_{j=1}^{[\sqrt{n}]}
		\dfrac{\delta_{j}(y)\cos j(t_k-y)}{(q^j+q^{-j})(1+\delta_{j}(y))}
\right)
\mathop{\text{sign}} \sin(ny-\frac{\beta\pi}{2})
=
\end{equation*}
\begin{equation}\label{s1}
=\left(
	\mathcal{P}_q(t_k-y)-\frac{1}{2}
\right)
\mathop{\text{sign}} \sin(ny-\frac{\beta\pi}{2})	
+\gamma_4(y)+\gamma_5(y),
\end{equation}
Із \eqref{SP_Psi}, \eqref{summa_gamma} та \eqref{s1}  отримуємо \eqref{SP_Phi_}.
Лему доведено.

\textbf{4. Доведення теореми 1.}
Відповідно до теореми 4 для доведення \eqref{dno1} і \eqref{dno2} достатньо показати, що для довільних $q\in(0,1)$, $\beta\in \mathbb{R}$ і всіх номерів $n\geqslant n_{q}$ ядра Неймана $N_{q,\beta}(t)$ задовольняють умову $C_{y_0,2n}$, де $y_0$ --- точка, в якій функція $|\Phi_{q,\beta,n}(\cdot)|$, де $\Phi_{q,\beta,n}(\cdot)=(N_{q,\beta}\ast\varphi_n)(\cdot)$, а $\varphi_n(\cdot)$ задана рівністю \eqref{varp}, досягає найбільшого значення, тобто
\begin{equation*}
|\Phi_{q,\beta,n}(y_0)|
	=
		|(N_{q,\beta}\ast\varphi_n)(y_0)|
	=
		\|N_{q,\beta}\ast\varphi_n\|_C.
\end{equation*}

Функція
\begin{equation*}
\Phi_{q,\beta,n}(t)=(N_{q,\beta}\ast\varphi_n)(t)=\dfrac{4}{\pi}\sum\limits_{\nu=0}^{\infty}\dfrac{q^{(2\nu+1)n}}{n(2\nu+1)^2}\sin\left((2\nu+1)nt-\dfrac{\beta\pi}{2}\right),
\end{equation*}
періодична з періодом $2\pi/n$ і така, що
$\Phi_{q,\beta,n}(t+\dfrac{\pi}{n})=-\Phi_{q,\beta,n}(t)$. Тому максимальне значення $\pi/n$-періодичної функції $|\Phi_{q,\beta,n}(\cdot)|$ на $[0,\dfrac{\pi}{n})$ досягається у точці $y_0=y_0(n,q,\beta)=\dfrac{\theta_n\pi}{n}$, де $\theta_n$ --- корінь рівняння \eqref{theta}, $\theta_n\in[0,1)$. Зауважимо, що цей корінь єдиний на $[0,1)$. Дійсно, розглянемо функції
\begin{equation}\label{Gq}
G_{q}(x)=
\sum\limits_{\nu=0}^{\infty}
	\dfrac{q^{(2\nu+1)n}}{(2\nu+1)n}\cos(2\nu+1)x
=\frac{1}{4}
	\ln\dfrac{1+2q^n\cos x+q^{2n}}{1-2q^n\cos x+q^{2n}},
\end{equation}
\begin{equation}\label{Hq}
H_{q}(x)=
\sum\limits_{\nu=0}^{\infty}
	\dfrac{q^{(2\nu+1)n}}{(2\nu+1)n}\sin(2\nu+1)x
=\frac{1}{2}
	\mathop{\mathrm{arctg}}\dfrac{2q^n\sin x}{1-q^{2n}}.
\end{equation}

З \eqref{Gq} випливає, що функція $G_{q}(x)$ спадає на $(0,\pi)$, зростає на $(\pi,2\pi)$ та 
\begin{equation*}
G_{q}(x)>0,\; x\in(0,\dfrac{\pi}{2})\cup(\dfrac{3\pi}{2},2\pi),
\end{equation*}
\begin{equation*}
G_{q}(x)<0,\; x\in(\dfrac{\pi}{2},\dfrac{3\pi}{2}).
\end{equation*}
В силу \eqref{Hq} функція $H_{q}(x)$ додатна на $(0,\pi)$ і від’ємна на $(\pi,2\pi)$.
Оскільки
\begin{equation*}
\sum\limits_{\nu=0}^{\infty}
	\dfrac
		{q^{(2\nu+1)n}}
		{(2\nu+1)n}
	\cos\left(
		(2\nu+1)\theta_n\pi-\dfrac{\beta\pi}{2}
	\right)
	=
		G_{q}(\theta_n\pi)\cos\dfrac{\beta\pi}{2}
		+
		H_{q}(\theta_n\pi)\sin\dfrac{\beta\pi}{2},
\end{equation*}
то, враховуючи зазначені вище властивості функцій $G_{q}(x)$ та $H_{q}(x)$ та рівність \eqref{theta}, отримуємо єдиність кореня $\theta_n$ на $[0,1)$ та наступні включення:
\begin{equation}\label{ny_0_1}
ny_0\in[\dfrac{\pi}{2},\pi) \text{ при } \beta\in[0,1)\cup[2,3),
\end{equation}
\begin{equation}\label{ny_0_2}
ny_0\in[0, \dfrac{\pi}{2}) \text{ при } {\beta\in[1,2)\cup[3,4)}.
\end{equation}

Згідно з рівностями \eqref{ny_0_1}, \eqref{ny_0_2}, рівністю (19) роботи \cite{Stepanets_1994} і лемою~2 з \cite{Stepanets_1994} для ядер $\Psi_{\beta}(t)=N_{q,\beta}(t)$ виконується умова $|\lambda_{j}(y_0)|\not=0$, ${j=\overline{1,n}}$. 
Тому для фундаментального $SK$-сплайна $\overline{S\Psi}_{\beta,1}(y,t)=\overline{SN}_{q,\beta,1}(y,t)$, породженого ядром Неймана $N_{q,\beta}(t)$, має місце представлення \eqref{SP_Phi_}.
Наступне твердження містить оцінку зверху суми $\sum\limits_{l=1}^5|\gamma_l(y_0)|$ для ядер $N_{q,\beta}(t)$.

\textbf{Лема 3.} \textit{Нехай величини $\gamma_l(y_0)$, $l=\overline{1,5}$, задаються рівностями \eqref{gamma_2}, \eqref{gamma_3}, \eqref{gamma_1}--\eqref{gamma_5}, в яких $\psi(n)=\dfrac{q^n}{n}$, $q\in(0,1)$, $\beta\in\mathbb{R}$. Тоді при $n\geqslant2$ та виконанні умов \eqref{umova_z} 
справедлива оцінка}
\begin{equation*}
\sum\limits_{l=1}^5|\gamma_l(y_0)|
\leqslant
\dfrac{24}{5(1-q)}q^{\sqrt{n}}
+\dfrac{160}{63}\dfrac{2\sqrt{n}-1}{n(\sqrt{n}-1)}\; \dfrac{q}{(1-q)^2}.
\end{equation*}

\textbf{Доведення.} Для оцінки кожного з доданків $|\gamma_l(y_0)|$, $l=\overline{1,5}$, нам будуть потрібні оцінки зверху величин $|r_{j}(y_0)|$ та $|R_{j}(y_0)|$ при $j=\overline{0,n-1}$. Знайдемо їх. З \eqref{r_nj_1} маємо
\begin{equation*}
|r_{j}^{(1)}(y_0)|\leqslant
\dfrac{q^{3n-j}}{(3n-j)^2}+\sum\limits_{m=2}^{\infty}\left(\dfrac{q^{(2m+1)n-j}}{((2m+1)n-j)^2}+\dfrac{q^{(2m-1)n+j}}{((2m-1)n+j)^2}\right)=
\end{equation*}
\begin{equation}\label{add_r_n}
=\sum\limits_{m=1}^{\infty}\left(\dfrac{q^{(2m+1)n-j}}{((2m+1)n-j)^2}+\dfrac{q^{(2m+1)n+j}}{((2m+1)n+j)^2}\right).
\end{equation}
Оскільки послідовність $\dfrac{q^k}{k^2}$ опукла, то виконується нерівність ${\dfrac{q^{k-j}}{(k-j)^2}+\dfrac{q^{k+j}}{(k+j)^2}<\dfrac{q^{k-n}}{(k-n)^2}+\dfrac{q^{k+n}}{(k+n)^2}}$, $k>n$, $j=\overline{0,n-1}$. Тому з \eqref{add_r_n} знаходимо
\begin{equation*}
|r_{j}^{(1)}(y_0)|
	\leqslant
		\sum\limits_{m=1}^{\infty}
			\left(
				\dfrac{q^{2mn}}{(2mn)^2}
				+
				\dfrac{q^{2(m+1)n}}{(2(m+1)n)^2}
			\right)
	=
\end{equation*}
\begin{equation}\label{|r_nj1|}
= \dfrac{q^{2n}}{4n^2}+\sum\limits_{m=2}^{\infty}\dfrac{q^{2mn}}{2m^2n^2}\leqslant
\dfrac{1}{4n^2}\sum\limits_{m=1}^{\infty}q^{2mn}=\dfrac{q^{2n}}{4n^2(1-q^{2n})}.
\end{equation}

З \eqref{cos_y}  та \eqref{r_nj_2} маємо
\begin{equation*}
|r_{j}^{(2)}(y_0)|
	\leqslant|
		\cos(ny_0-\frac{\beta\pi}{2})|
		\left(
			\dfrac{q^{n-j}}{(n-j)^2}
			-\dfrac{q^{n+j}}{(n+j)^2}
		\right)
	\leqslant
\end{equation*}
\begin{equation}\label{|r_nj2|}
	\leqslant
		\dfrac{q^{2n}}{3(1-q^{2n})}
		\left(
			q
			-\dfrac{q^{2n-1}}{(2n-1)^2}
		\right).
\end{equation}

Із  \eqref{a_n} та \eqref{r_nj_3} знаходимо
\begin{equation}\label{|alpha_n|}
|r_{j}^{(3)}(y_0)|
	\leqslant
		\dfrac{q^{2n}}{3(1-q^{2n})}
		\left(
			q
			+\dfrac{q^{2n-1}}{(2n-1)^2}
		\right).
\end{equation}

Об’єднавши \eqref{|r_nj1|}, \eqref{|r_nj2|} та \eqref{|alpha_n|}, для величини $r_{j}(y_0)$ отримуємо оцінку
\begin{equation}\label{|r_nj|}
|r_{j}(y_0)|
	\leqslant 
		|\sum\limits_{\nu=1}^{3}
			r_{j}^{(\nu)}(y_0)|
	\leqslant 
		\dfrac{q^{2n}}{1-q^{2n}}
		\left(
			\frac{2q}{3}
			+\dfrac{1}{4n^2}
		\right)
	\leqslant 
		\frac{3}{4}\dfrac{q^{2n}}{1-q^{2n}},
j=\overline{0,n-1}.
\end{equation}

При $j=0$ оцінку \eqref{|r_nj|} можна покращити. Дійсно, в силу \eqref{r_n1} маємо
\begin{equation*}
|r_{0}^{(1)}(y_0)|\leqslant2\sum\limits_{m=2}^{\infty}\dfrac{q^{(2m-1)n}}{((2m-1)n)^2}\leqslant
\dfrac{2}{9n^2}\sum\limits_{m=2}^{\infty}q^{(2m-1)n}=\dfrac{2}{9n^2}\dfrac{q^{3n}}{1-q^{2n}},
\end{equation*}
а з \eqref{r_n3} та \eqref{a_n}
\begin{equation*}
\left| r_{0}^{(3)}(y_0)  \right|\leqslant \dfrac{2}{3n^2} \dfrac{q^{3n}}{1-q^{2n}}.
\end{equation*}
Тоді, враховуючи \eqref{r_n2},
\begin{equation}\label{|r_n|}
|r_{0}(y_0)|\leqslant |r_{0}^{(1)}(y_0)+r_{0}^{(3)}(y_0)| \leqslant \dfrac{8}{9n^2} \dfrac{q^{3n}}{1-q^{2n}}.
\end{equation}

Із \eqref{lambda_n-j} для величини $|\lambda_{n-j}(y_0)|$ отримуємо зображення
\begin{equation*}
|\lambda_{n-j}(y_0)|
	=
		\left|
			\mathop{\text{sign}} \sin(ny-\frac{\beta\pi}{2})
			\left(
				\dfrac
					{q^{n-j}}
					{(n-j)^2}
				+
				\dfrac
					{q^{n+j}}
					{(n+j)^2}
			\right)
			+
			r_{j}(y_0)
		\right|,
\end{equation*}
з якого безпосередньо випливає оцінка
\begin{equation}\label{|lambda_n-j|0'}
|\lambda_{n-j}(y_0)|\leqslant \dfrac{q^{n-j}}{(n-j)^2}+\dfrac{q^{n+j}}{(n+j)^2}+|r_{j}(y_0)|.
\end{equation}
Оскільки внаслідок \eqref{ny_0_1} і \eqref{ny_0_2}
\begin{equation}\label{sign_not_zero}
\sin(ny_0-\dfrac{\beta \pi}{2})=\sin ny_0\cos\dfrac{\beta \pi}{2}-\cos ny_0\sin\dfrac{\beta \pi}{2}\not=0.
\end{equation}
то отримуємо також оцінку
\begin{equation}\label{|lambda_n-j|0''}
|\lambda_{n-j}(y_0)|\geqslant \dfrac{q^{n-j}}{(n-j)^2}+\dfrac{q^{n+j}}{(n+j)^2}-|r_{j}(y_0)|.
\end{equation}

В силу \eqref{|lambda_n-j|}, \eqref{|lambda_n-j|0'} та \eqref{|lambda_n-j|0''}
\begin{equation}\label{|R_n,j|}
|R_{j}(y_0)|\leqslant |r_{j}(y_0)|,\; j=\overline{0,n-1}.
\end{equation}

Перейдемо до оцінки величини $|\gamma_1(y_0)|$.
Взявши до уваги оцінки  \eqref{|lambda_n-j|0''} та \eqref{|r_nj|}, маємо
\begin{equation*}
|\lambda_{n-j}(y_0)|
	\geqslant
\end{equation*}
\begin{equation*}
	\geqslant 	
		\dfrac{q^{n-j}}{(n-j)^2}
		+\dfrac{q^{n+j}}{(n+j)^2}
		-\frac{3}{4}\dfrac{q^{2n}}{1-q^{2n}}
	=
\end{equation*}
\begin{equation}\label{|lambda_n-j|3}
=
	\dfrac{q^{n}}{(n-j)^2} 
	\left(
		q^{-j}
		+\dfrac{(n-j)^2}{(n+j)^2}q^{j}
		-\frac{3(n-j)^2}{4}\dfrac{q^{n}}{1-q^{2n}}
	\right).
\end{equation}

Оскільки при $j=\overline{0,n-1}\;$ 
$\dfrac{2}{15(n-j)^2}q^{-j}
>\dfrac{2}{15n^2}
>\dfrac{2q^{\sqrt{n}}}{15n^2}$, то з умови \eqref{umova_z} випливає нерівність
\begin{equation*}
\dfrac{2}{15(n-j)^2}q^{-j}
	>\dfrac{q^{n}}{1-q^{2n}},
\end{equation*}
яка еквівалентна наступній нерівності:
\begin{equation}\label{um1}
\dfrac{q^{-j}}{10}
	>\frac{3(n-j)^2}{4}\dfrac{q^{n}}{1-q^{2n}},
\; j=\overline{0,n-1}.
\end{equation}

В силу \eqref{um1} виконуються оцінки
\begin{equation*}
q^{-j}
+\dfrac{(n-j)^2}{(n+j)^2}q^{j}
-\frac{3(n-j)^2}{4}\dfrac{q^{n}}{1-q^{2n}}
=
\end{equation*}
\begin{equation*}
=
 \dfrac{9q^{-j}}{10}
 +\dfrac{q^{-j}}{10}
 +\dfrac{(n-j)^2}{(n+j)^2}q^{j}
 -\frac{3(n-j)^2}{4}\dfrac{q^{n}}{1-q^{2n}}
>
\end{equation*}
\begin{equation}\label{ots}
>
	\dfrac{9q^{-j}}{10},
\; j=\overline{0,n-1}.
\end{equation}

Об’єднуючи \eqref{|lambda_n-j|3} та \eqref{ots}, маємо
\begin{equation}\label{mod_lambda}
|\lambda_{n-j}(y_0)|
	\geqslant
		\dfrac{9q^{n-j}}{10(n-j)^2}.
\end{equation}

Враховуючи, що для $x\in[0,\dfrac \pi2)$ справджується нерівність $\cos x\geqslant 1-\dfrac {2x}{\pi}>0$, отримуємо
\begin{equation}\label{cos1}
\cos \dfrac{j\pi}{2n}\geqslant 1-\dfrac {j}{n}=\dfrac {n-j}{n},\; j=\overline{0,n-1}.
\end{equation}

З \eqref{mod_lambda} та \eqref{cos1} маємо
\begin{equation}\label{|lambda_n-j|_2}
\dfrac{n^2}{q^{n}}|\lambda_{n-j}(y_0)|^2\cos\dfrac{j\pi}{2n} > \dfrac{81n}{100(n-j)^3}q^{n-2j}.
\end{equation}

З \eqref{z_nj} і \eqref{|R_n,j|} випливає, що $|z_{j}(y_0)|\leqslant 2|r_{j}(y_0)|$. Тому враховуючи \eqref{|r_nj|}, \eqref{|lambda_n-j|_2} та умову \eqref{umova_z}, з \eqref{gamma_2} одержуємо
\begin{equation*}
\left|\gamma_1(y_0)\right|
	\leqslant
		\dfrac{400}{81}\max_{0\leqslant j\leqslant n-1}| r_{j}(y_0)| \dfrac{1}{q^{n}}
		\sum_{j=0}^{n-1}
			\dfrac{(n-j)^3}{n}q^{2j}
	<
		\dfrac{100\,n^2q^{n}}{27(1-q^{2n})}
		\sum_{j=0}^{\infty}
			q^{2j}
	\leqslant
\end{equation*}
\begin{equation}\label{r1}
\leqslant
\dfrac{40}{81}q^{\sqrt{n}}\;\dfrac{1}{1-q^2}.
\end{equation}

Оцінимо $|\gamma_2(y_0)|$. З умови \eqref{umova_z} при $n\geqslant2$ випливає нерівність
\begin{equation*}
q^{2n}\leq \dfrac{1}{900}.
\end{equation*}

Тоді з \eqref{gamma_3}, \eqref{|R_n,j|}, \eqref{|r_n|} і \eqref{umova_z}  отримуємо
\begin{equation*}
|\gamma_2(y_0)|
	\leqslant
		\dfrac
		{
			 \dfrac{8q^{2n}}{9(1-q^{2n})}
		}
		{
			2\left|2-\dfrac{8q^{2n}}{9(1-q^{2n})}\right|
		}
	=
		\dfrac{2q^{2n}}{9-13q^{2n}}
	=
		\dfrac{1-q^{2n}}{9-13q^{2n}}\;\dfrac{2q^{2n}}{1-q^{2n}}=
\end{equation*}
\begin{equation*}
	=
		\left(
			\dfrac {1}{13}
			+
			\dfrac{4}{13(9-13q^{2n})}
		\right)
		\dfrac{2q^{2n}}{1-q^{2n}}
	<
\end{equation*}
\begin{equation}\label{r2_2}
	<
		\dfrac {3}{26} \dfrac{2q^{2n}}{1-q^{2n}}
	<
		\dfrac{2q^{n+\sqrt{n}}}{65n^2}.
\end{equation}

Оцінимо величину $|\gamma_3(y_0)|$.
Взявши до уваги оцінки \eqref{|r_nj|}, \eqref{|lambda_n-j|0''}, \eqref{cos1} та \eqref{ots}, маємо
\begin{equation*}
\dfrac{n^2}{q^{n}}|\lambda_{n-j}(y_0)|\cos\dfrac{j\pi}{2n}
	\geqslant 
		\dfrac{n^2}{q^{n}} 
		\left(
			\dfrac{q^{n-j}}{(n-j)^2}
			+
			\dfrac{q^{n+j}}{(n+j)^2}
			-
			\frac{3}{4}\dfrac{q^{2n}}{1-q^{2n}}
		\right)
		\dfrac {n-j}{n}
	=
\end{equation*}
\begin{equation}\label{|lambda_n-j|2}
	=
		\dfrac{nq^{-j}}{n-j}
		+
		\dfrac{n(n-j)}{(n+j)^2}q^{j}
		-
		\frac{3n(n-j)}{4}\dfrac{q^{n}}{1-q^{2n}}
	> 
		\dfrac{9nq^{-j}}{10(n-j)}.
\end{equation}

Тому, взявши до уваги \eqref{|lambda_n-j|2}, з \eqref{gamma_1} знаходимо
\begin{equation*}
\left|\gamma_3(y_0)\right|
	<
		\dfrac{20}{9}
		\sum_{j=[\sqrt{n}]+1}^{n-1}
			\dfrac{n-j}{n}q^j
	<
\end{equation*}
\begin{equation}\label{r3}
	<
		\dfrac{20}{9}
		\sum_{j=[\sqrt{n}]+1}^{n-1}
		q^j
	=
		\dfrac{20(q^{[\sqrt{n}]+1}-q^{n-1})}{9(1-q)}
	\leqslant
		\dfrac{20q^{\sqrt{n}}}{9(1-q)}.
\end{equation}

Щоб оцінити величину $|\gamma_4(y_0)|$ спочатку оцінимо зверху величину $|\delta_j(y_0)|$ вигляду \eqref{delta_0}. В силу \eqref{|lambda_n-j|}
\begin{equation*}
\dfrac{n^2}{q^{n}}
	|\lambda_{n-j}(y_0)|
		\cos\dfrac{j\pi}{2n}
	=
\end{equation*}
\begin{equation*}
	=
		\left(
			\dfrac{n^2}{(n-j)^2}\dfrac{q^{n-j}}{q^{n}}
			+
			\dfrac{n^2}{(n+j)^2}\dfrac{q^{n+j}}{q^{n}}
			+
			R_{j}(y_0)\dfrac{n^2}{q^{n}}
		\right)
		\cos\dfrac{j\pi}{2n}
	=
\end{equation*}
\begin{equation*}
	=
		\left(
			(1+\dfrac{j(2n-j)}{(n-j)^2})q^{-j}
			+
			(1-\dfrac{j(2n+j)}{(n+j)^2})q^{j}
			+
			R_{j}(y_0)\dfrac{n^2}{q^{n}}
		\right)
		\cos\dfrac{j\pi}{2n}
	=
\end{equation*}
\begin{equation*}
	=
		(q^{-j}+q^{j})(1-2\sin^2\dfrac{j\pi}{4n})
	+
\end{equation*}
\begin{equation}\label{cos_lambda0}
	+
		\left(
			\dfrac{j(2n-j)}{(n-j)^2}q^{-j}
			-
			\dfrac{j(2n+j)}{(n+j)^2}q^{j}
			+
			R_{j}(y_0)\dfrac{n^2}{q^{n}}
		\right)
		\cos\dfrac{j\pi}{2n}.
\end{equation}

З \eqref{delta_0}, \eqref{|r_nj|}, \eqref{|R_n,j|}, \eqref{cos_lambda0} із врахуванням опуклості послідовності $q^k$ для величин $|\delta_{j}(y_0)|$ випливають нерівності
\begin{equation*}
|\delta_{j}(y_0)|
	\leqslant 
		2\sin^2\dfrac{j\pi}{4n}
		+
		\dfrac{1}{q^{-j}+q^{j}}
		\left(
			\dfrac{j(2n-j)}{(n-j)^2}q^{-j}
			+
			\dfrac{j(2n-j)}{(n-j)^2}q^{j}
			+
			|R_{j}(y_0)|\dfrac{n^2}{q^{n}}
		\right)
	\leqslant
\end{equation*}
\begin{equation*}
	\leqslant 
		2\left(
			\dfrac{j\pi}{4n}
		\right)^2
		+
		\dfrac{j(2n-j)}{(n-j)^2}
		+
		\dfrac{n^2|r_{j}(y_0)|}{q^{n-j}+q^{n+j}}
	\leqslant 
		\dfrac{j^2\pi^2}{8n^2}
		+
		\dfrac{j(2n-j)}{(n-j)^2}
		+
		\frac{3n^2}{8}\dfrac{q^{n}}{1-q^{2n}}
	=
\end{equation*}
\begin{equation}\label{delta0}
	= 
		\dfrac{8j(2n-j)}{7(n-j)^2}
		+
		\left(
			\dfrac{j^2\pi^2}{8n^2}
			+
			\frac{3n^2}{8}\dfrac{q^{n}}{1-q^{2n}}
			-
			\dfrac{j(2n-j)}{7(n-j)^2}
		\right).
\end{equation}

При кожному фіксованому $n$ функція $f_n(j)=\dfrac{j(2n-j)}{7(n-j)^2}-\dfrac{j^2\pi^2}{8n^2}$ зростає. 
Дійсно
\begin{equation*}
\frac{d}{dj}
\left(
	\dfrac
		{j(2n-j)}
		{7(n-j)^2}
	-
	\dfrac
		{j^2\pi^2}
		{8n^2}
\right)
	=
		\dfrac
			{2n^{2}}
			{7(n-j)^{3}}
		-
		\dfrac
			{j\pi^2}
			{4n^2}
	=
\end{equation*}
\begin{equation*}
	=
		\dfrac
			{8n^{4}-7j\pi^2(n-j)^{3}}
			{28n^2(n-j)^{3}}
	>0
\end{equation*}
(оскільки функція $g_n(x)=8n^{4}-7x\pi^2(n-x)^{3}$ в точці мінімуму $x=\frac{1}{4}n$ набуває додатне значення).
Тому, з врахуванням \eqref{umova_z}, при $j=\overline{1,[\sqrt{n}]}$ та $n\geqslant2$ отримуємо
\begin{equation*}
\dfrac{j(2n-j)}{7(n-j)^2}
-
\dfrac{j^2\pi^2}{8n^2}
	\geqslant
		\dfrac{2n-1}{7(n-1)^2}
		-
		\dfrac{\pi^2}{8n^2}
	>
		\frac{3n^2}{8}\dfrac{q^{n}}{1-q^{2n}}.
\end{equation*}

Отже, вираз у дужках в правій частині \eqref{delta0} від’ємний. Тоді з \eqref{delta0} випливає, що
\begin{equation}\label{delta}
|\delta_{j}(y_0)|
	\leqslant 
		\dfrac{8j(2n-j)}{7(n-j)^2}.
\end{equation}

Формули \eqref{gamma_4}, \eqref{|lambda_n-j|2} та \eqref{delta} дозволяють одержати 
при $n\geqslant2$ наступну оцінку величини $\gamma_4(y_0)$:
\begin{equation*}
|\gamma_4(y_0)|
	\leqslant
		2\sum_{j=1}^{[\sqrt{n}]}
			\dfrac
			{
				\dfrac{8j(2n-j)}{7(n-j)^2}
			}
			{
				\dfrac{9nq^{-j}}{10(n-j)}
			}
	=
		\frac{160}{63n}
		\sum_{j=1}^{[\sqrt{n}]}
			\dfrac{j(2n-j)}{n-j}q^j
	\leqslant
\end{equation*}
\begin{equation*}
	\leqslant
		\dfrac{160}{63}
		\dfrac{2n-\sqrt{n}}{n(n-\sqrt{n})}
		\sum_{j=1}^{[\sqrt{n}]}
			jq^j
	<
		\dfrac{160}{63}
		\dfrac{2\sqrt{n}-1}{n(\sqrt{n}-1)}
		\sum_{j=1}^{\infty}
			jq^j
	<
\end{equation*}
\begin{equation}\label{|R|}
	<
		\dfrac{160}{63}
		\dfrac{2\sqrt{n}-1}{n(\sqrt{n}-1)}
		\; \dfrac{q}{(1-q)^2}.
\end{equation}

В силу \eqref{gamma_5} для величини $|\gamma_5(y_0)|$ маємо
\begin{equation}\label{r4}
\left|\gamma_5(y_0)\right|
	\leqslant 
		2\sum\limits_{j=[\sqrt{n}]+1}^{\infty}
			q^j
	=
		2\dfrac{q^{[\sqrt{n}]+1}}{1-q}
	<
		2\dfrac{q^{\sqrt{n}}}{1-q}.
\end{equation}

Об’єднавши оцінки \eqref{r1}, \eqref{r2_2}, \eqref{r3}, \eqref{|R|} та \eqref{r4}, при $n\geqslant2$ одержимо
\begin{equation*}
\sum_{k=1}^5
	|\gamma_k(y_0)|
	<
\end{equation*}
\begin{equation*}
	<
		\dfrac{40}{81}q^{\sqrt{n}}
		\;\dfrac{1}{1-q^2}
		+
		\dfrac{2q^{n+\sqrt{n}}}{65n^2}
		+
		\dfrac{20q^{\sqrt{n}}}{9(1-q)}
		+
		\dfrac{160}{63}
		\dfrac{2\sqrt{n}-1}{n(\sqrt{n}-1)}
		\; \dfrac{q}{(1-q)^2}
		+
		\dfrac{2q^{\sqrt{n}}}{1-q}
	<
\end{equation*}
\begin{equation*}
	<
		\dfrac{q^{\sqrt{n}}}{1-q}
		(
			0{,}494
			+
			0{,}008
			+
			2{,}223
			+
			2
		)
		+
		\dfrac{160}{63}
		\dfrac{2\sqrt{n}-1}{n(\sqrt{n}-1)}
		\; \dfrac{q}{(1-q)^2}
	<
\end{equation*}
\begin{equation*}
	<
		\dfrac{24}{5(1-q)}q^{\sqrt{n}}
		+
		\dfrac{160}{63}
		\dfrac{2\sqrt{n}-1}{n(\sqrt{n}-1)}
		\; \dfrac{q}{(1-q)^2}.
\end{equation*}

Лему доведено.

Згідно з лемою~2 роботи \cite{My_JAT} для довільного $x\in\mathbb{R}$ і довільного $q\in(0,1)$
\begin{equation}\label{f_x}
\mathcal{P}_q(x)
	>
		\left(
			\dfrac{1}{2}
			+
			\dfrac{2q}{(1+q^2)(1-q)}
		\right)
		\left(
			\dfrac{1-q}{1+q}
		\right)^{\frac {4}{1-q^2}}.
\end{equation}
Тому з леми 3 та нерівності \eqref{f_x} випливає, що при $n\geqslant2$ за умов \eqref{umova_n_0} та \eqref{umova_z}
\begin{equation}\label{geq0}
\mathcal{P}_q(t_k-y_0)
+
\sum\limits_{m=1}^5
	\gamma_m(y_0)
\mathop{\text{sign}} \sin(ny_0-\frac{\beta\pi}{2})
	\geqslant
		0.
\end{equation}
В силу зображення \eqref{SP_Phi_}, а також \eqref{sign_not_zero} і нерівності \eqref{geq0} робимо висновок, що при $n\geqslant2$ за умов \eqref{umova_n_0} та \eqref{umova_z} справедливе включення $N_{q,\beta}(t)\in C_{y_0,2n}$.
Теорему доведено.

\renewcommand{\refname}{}
\makeatletter\renewcommand{\@biblabel}[1]{#1.}\makeatother

\end{document}